\documentclass[10pt]{article}
\usepackage{amsmath, amsfonts, amsthm, amssymb, verbatim}
\usepackage{graphicx}
\usepackage[mathcal]{euscript}
\usepackage{amstext}
\usepackage{amssymb,mathrsfs}
\usepackage[latin1]{inputenc}
\newcommand \be     {\begin{equation}}
\newcommand \ee     {\end{equation}}

\newcommand \del     \partial

\def\XXint#1#2#3{{\setbox0=\hbox{$#1{#2#3}{\int}$}
\vcenter{\hbox{$#2#3$}}\kern-.5\wd0}}

\DeclareMathOperator\dive {div}
 \newcommand{\R}{\mathbb R}

\newtheorem{theorem}{Theorem}[section]

 \newtheorem{remark}[theorem]{Remark}
\newtheorem{lemma}[theorem]{Lemma}

 \def\beqs{\begin{eqnarray*}}
 \def\enqs{\end{eqnarray*}}
 \def\beq{\begin{eqnarray}}
 \def\enq{\end{eqnarray}}

\begin{document}
\title{The Hardy and Caffarelli-Kohn-Nirenberg Inequalities Revisited}
\author{Aldo Bazan$^1$, Wladimir Neves$^1$}

\date{}

\maketitle

\footnotetext[1]{ Instituto de Matem\'atica, Universidade Federal
do Rio de Janeiro, C.P. 68530, Cidade Universit\'aria 21945-970,
Rio de Janeiro, Brazil. E-mail: {\sl aabp2003@pg.im.ufrj.br,
wladimir@im.ufrj.br.}
\newline
\textit{To appear in:}
\newline
\textit{Key words and phrases.} Gauss-Green Theorem, Functional
inequalities, singularities.}

%
%
\begin{abstract}
In this paper some important inequalities are revisited. First, as
motivation, we give another proof of the Hardy's inequality
applying convenient vector fields as introduced by Mitidieri, see
\cite{MTD}. Then, we investigate a particular case of the
Caffarelli-Kohn-Nirenberg's inequality. Finally, we study the
Rellic's inequality.
\end{abstract}
%

\maketitle

\tableofcontents

\section{Introduction} \label{IN}

We begin our study by Hardy's inequality, in fact as motivation.
Another proof of this inequality is given applying the original
ideas of convenient vector fields as introduced by Mitidieri, see
\cite{MTD}. Although, differently from that paper, we use during
the proof the generalized Young's inequality, which gives to us a
simple way (a la Calculus) to obtain the best constants in some
sort of this inequalities.

\smallskip
Then, we investigate a particular case of the
Caffarelli-Korn-Nirenberg's inequality. In fact, we prove this
inequality by an interpolation argument based on two convenient
parameter points, see Theorem \ref{CKNT}. That is, first we prove
Lemma \ref{L1}, concerning the Caffarelli-Korn-Nirenberg
inequality for $b= a+1$, applying the idea of convenient vector
fields, where the technic gives to us the best constant. After, we
show a second lemma for $b= a$, where the Sobolev's inequality is
used with a convenient function and, further we apply the result
proved in the previous Lemma \ref{L1}. To our knowledge this
procedure is completely new.

\smallskip
Finally, we investigate the Rellich's inequality, which is a
second order type inequality like a generalization of Hardy's one.
Our proof is based on the considered particular case of
Caffarelli-Korn-Nirenberg's inequality for $b= a+1$, with $a=1$.
Again, the best constant is recovered in our analysis due to
critical point procedure.

\subsection{Main purpose}

In this paper, we first consider the following sharp version of
the Hardy's inequality
$$
  \int_{\R^{n}}\frac{\left|u(x)\right|^{p}}{\left\|x\right\|^{p}}\, dx
  \leq \left(\frac{p}{n-p}\right)^{p}
  \int_{\R^{n}}\left\|\nabla u(x)\right\|^{p} \, dx,
$$
where
$u$ is a $C^{\infty}_{c}(\R^{n})$ function and $p$ is a real
number, such that $1 < p < n$.
Moreover, for $u \in C_c(\R^n \setminus \{0\})$ and $1 < n < p$,
the sharp version of Hardy's turns
$$
  \int_{\R^{n}}\frac{\left|u(x)\right|^{p}}{\left\|x\right\|^{p}}\, dx
  \leq \left(\frac{p}{p-n}\right)^{p}
  \int_{\R^{n}}\left\|\nabla u(x)\right\|^{p} \, dx.
$$
In one dimensional case, we show the following version of Hardy's
inequality
$$
  \int_{0}^{\infty} \left(\frac{\eta(x)}{x}\right)^{p} \, dx
  <
  \left(\frac{p}{p-1}\right)^{p}\int_{0}^{\infty} u(x)^{p} \, dx,
$$
where $p > 1$, $u \in C_{c}^{\infty}(\R+)$ is a nonnegative and
non-identically zero function, and
$$
  \eta(x):= \int_{0}^{x} u(t) \, dt.
$$
The reader is addressed to Section \ref{FN} for the functional
notation. Following the original idea of Mitidiere \cite{MTD}, we
give a new proof of these inequalities applying a nice and simple
technique, which further possibility us to recovery direct the
best constant.

\medskip
Also, we analyze a particular case of the well-known inequality
due to Caffarelli-Kohn-Nirenberg, which asserts that
$$
  \left(\int_{\R^{n}}\frac{\left|u(x)\right|^{p}}{\left\|x\right\|^{bp}}dx\right)^{2/p}
  \leq C(n,b,p) \int_{\R^{n}}\frac{\left\|\nabla
  u(x)\right\|^{2}}{\left\|x\right\|^{2a}}dx,
$$
where $C(n,b,p)$ is a positive constant, $n\geq 3$, $u\in
C^{\infty}_{c}(\R^{n})$ and also
$$-\infty<a<\frac{n-2}{2}, \quad
a\leq b\leq a+1, \quad p=\frac{2n}{n-2+2(b-a)}.
$$
Following the same ideas applied to prove the Hardy inequality,
first we were able to show the Caffarelli-Kohn-Nirenberg's
inequality considered, when $b=a+1$. In this particular case, the
sign of $a$ has no influence and we recover the sharp constant,
that is
$$
  C(n,b,2)= \frac{4}{(n-2-2a)^{2}}.
$$
On the other hand, when $b=a$ the sign of $a$ has to be considered
and we have differences between the constants, see \cite{CATRINA}.
Indeed, for $a<0$ the sharp constant is the well known sharp
constant of the Sobolev's inequality, that is
$$
  C(n,b,2^*)= K(n,2)^2,
$$
and when $a>0$ the constant is
$$
  C(n,b,2^*)= K(n,2)^2 \, \Big( 1 + a \sqrt{C(n,b,2)} \Big)^{2- 2/n}.
$$
The general case of Caffarelli-Kohn-Nirenberg'inequality
considered is proved by an interpolation argument between these
two previous cases, i.e. $b=a+1$ and $b=a$.

\medskip
At the end of the paper, we prove the Rellich's inequality, which
is a first order generalization of Hardy's inequality when $p=2$,
that is
$$
\label{RR}
    \int_{\R^{n}}\frac{\left|u(x)\right|^{2}}{\left\|x\right\|^{4}} \,
    dx
    \leq
    \left(\frac{4}{n(n-4)}\right)^2 \int_{\R^{n}}\left|\Delta
u(x)\right|^{2}dx,
$$
where $u \in C_c(\R^n)$, $n>4$ and
$\left(\frac{4}{n(n-4)}\right)^2$ is the sharp constant.

\medskip
\begin{remark}
%
%
Considering $p=2$, the inequality due to Hardy present here is a
particular case of the Caffarelli-Kohn-Nirenberg's inequality.
Indeed, it is enough to take $b=1$ and $a=0$  (and consequently
$p=2$).
\end{remark}

\medskip
An outline of this paper follows. In the rest of this section we
fix some definitions and notation. Moreover, we recall some
well-known results. The Hardy's inequality is proved on Section
\ref{HIS}. In Section \ref{CKNI}, first we prove two Lemmas, which
are the Caffarelli-Khon-Niremberg's inequality for $b=a+1$ and
$b=a$ respectively. Then, applying an interpolation argument we
were able to prove the general case. Finally, we show in Section
\ref{RIS} the Rellich's inequality.

\subsection{Functional notation and background} \label{FN}

At this point we fix the functional notation used in the paper and
recall some well known results.

\medskip
By $dx$ we denote the Lebesgue measure on $\R^n$. Moreover, we
denote by $\| . \|$ and $| . |$ respectively the Euclidean norm in
$\R^{n}$ and the absolute value in $\R$.

\medskip
We recall the generalized Young's inequality: For $\lambda
>0$ and all $V, W \in \R^n$, we have
\begin{equation}
\label{YI}
  V \cdot W \leq \lambda^{-p} \frac{\|V\|^p}{p} + \lambda^q
  \frac{\|W\|^q}{q},
\end{equation}
where $p,q \geq 1$ satisfying
$$
  \frac{1}{p} + \frac{1}{q}= 1, \quad \text{\big(or $q= p \, (q-1)$, or $p= q \, (p-1)\big)$}.
$$

\medskip
For any $U \subset \R^n$ the set $C_c^{\infty}(U)$ stands for the
space of all $C^\infty$ functions on $\R^n$ whose support is
compact and contained in $U$.
%
The Sobolev space $W^{1,p}(\R^n)$ is the set of all functions in
$L^p(\R^n)$ with first derivatives also in $L^p(\R^n)$, $(1 \leq p
< \infty)$, where the derivatives should be understood in the
usual weak sense.

\medskip
For $1 \leq p < n$, we set
$$
  p^*:= \frac{n p}{n-p},
$$
called the Sobolev conjugate of $p$. Thus the Sobolev's inequality
asserts that, for all functions $f \in W^{1,p}(\R^n)$
\begin{equation}
\label{SI}
   \left(\int_{\R^{n}}|f(x)|^{p^*}dx\right)^{1/p^*} \leq
   K(n,p)^{2} \, \left(\int_{\R^{n}}\left\|\nabla f(x)\right\|^{p} \,
   dx \right)^{1/p},
\end{equation}
where $K(n,p)$ is the sharp constant, given by
$$
\begin{aligned}
  K(n,p)&= \frac{1}{2^{1/n} \, \pi^{1/2} \, n}
             \left(\frac{p-1}{n-p}\right)^{1 - 1/p}
             \left(\frac{p}{p-1}\right)^{1/n} \\[5pt]
             &\times
             \left( \frac{\Gamma(n/2) \, \Gamma (n)}
             {\Gamma(n/p) \, \Gamma(n (1-1/p))}\right)^{1/n},
\end{aligned}
$$
and $\Gamma(s)$ is the Gamma-function.

\section{The Hardy inequality} \label{HIS}

The proof of the
Hardy's inequality follows with a nice strategy, which allow us to
apply the Gauss-Green Theorem. Then, the Young's inequality is
used to obtain our result.

\bigskip
\subsection{The case $(p \neq n, n > 1)$}

\begin{theorem}
Let $u$ be a function in $C_c^{\infty}(\R^{n})$ and $1 < p < n$.
Then,
\begin{equation}
\label{HR}
\int_{\R^{n}}\frac{\left|u(x)\right|^{p}}{\left\|x\right\|^{p}}dx
\leq \left(\frac{p}{n-p}\right)^{p}\int_{\R^{n}}\left\|\nabla
u(x)\right\|^{p}dx,
\end{equation}
where $\left(\frac{p}{n-p}\right)^{p}$ is the sharp constant.
\end{theorem}


\begin{proof} 1. First, let $V:\R^{n} \setminus \left\{0\right\} \to \R^n$ be a smooth vector
field, defined by
\begin{equation}
\label{VV}
       V(x):=\frac{x}{(p-n)\left\|x\right\|^{p}}.
\end{equation}
For each $i,j =1, \ldots,n$, this function verifies
$$
  \frac{\partial V_i(x)}{\partial x_j}= \frac{1}{(p-n) \, \|x\|^p}
  \; \delta_{ij} \, - \frac{p}{(p-n) \, \|x\|^{p-2}}
  \; x_i \, x_k \; \delta_{jk},
$$
where the usual summation convention and Kronecker delta notation
is used. Consequently, we have
$$
  \dive V(x)= -\frac{1}{\left\|x\right\|^{p}}.
$$

2. Now, the integral on the left side of $(\ref{HR})$ can be
rewritten in the following way
$$
  \int_{\R^{n}}\frac{\left|u(x)\right|^{p}}{\left\|x\right\|^{p}} \, dx
  =-\int_{\R^{n}}\left|u(x)\right|^{p} \dive V(x) \, dx.
$$
Then applying the Gauss-Green Theorem, we obtain
$$
  \begin{aligned}
  \int_{\R^{n}}\frac{\left|u(x)\right|^{p}}{\left\|x\right\|^{p}} \,
  dx
  &=p \int_{\R^{n}}  \left(\left|u(x)\right|^{p-1}V(x)\right) \cdot \left(\nabla \left|u(x)\right|\right) \,
  dx \\[5pt]
  &\leq p \int_{\R^n} \frac{\lambda^q}{q} \, |u(x)|^{(p-1)q}
  \|V(x)\|^q \, dx +  \int_{\R^n} \lambda^{-p} \|\nabla u(x)\|^p \, dx,
\end{aligned}
$$
where we have used the Young's inequality \eqref{YI}. Therefore,
it follows that
\begin{equation}
\label{BCR}
  \begin{aligned}
  \int_{\R^{n}}\frac{\left|u(x)\right|^{p}}{\left\|x\right\|^{p}}\, dx - \frac{p \, \lambda^{q}}{q}
  \int_{\R^{n}}\left|u(x)\right|^{p} \, \left\|V(x)\right\|^{q} \, dx &\leq
  \lambda^{-p}\int_{\R^{n}}\left\|\nabla u(x)\right\|^{p}dx \\[5pt]
  \left( 1-\frac{\lambda^{q}p}{q(n-p)^{q}} \right)
  \int_{\R^{n}}\frac{\left|u(x)\right|^{p}}{\left\|x\right\|^{p}}dx
  &\leq
  \lambda^{-p}\int_{\R^{n}}\left\|\nabla u(x)\right\|^{p} \, dx.
  \end{aligned}
\end{equation}
From \eqref{BCR} and a simple algebraic manipulation, we obtain
$$
  \int_{\R^{n}}\frac{\left|u(x)\right|^{p}}{\left\|x\right\|^{p}}\, dx
  \leq f(\lambda;n,p,q) \int_{\R^{n}}\left\|\nabla u(x)\right\|^{p} \,
  dx,
$$
where
\begin{equation}
\label{Hf}
  f(\lambda;n,p,q):= \frac{q (n-p)^q}{\lambda^p \big(q (n-p)^q -
  \lambda^q p \big)}.
\end{equation}

3. Finally, we proceed to obtain the sharp constant. Fixed $n$,
$p$ and thus $q$, we set the positive constant $\kappa= q
(n-p)^q$. Then, we have
$$
  f(\lambda)= \frac{\kappa}{\lambda^p \, \kappa - \lambda^{p+q} \, p}.
$$
So, we can derive $f$ and make it equal zero, to obtain a minimal
point candidate, that is
$$
  \lambda_0= \Big( \frac{\kappa}{p+q} \Big)^{1/q}.
$$
In fact, a straightforward calculation shows that $\lambda_0$ is
the point of minimum and
$$
  f(\lambda_0)= \left(\frac{p}{n-p}\right)^{p},
$$
which is the sharp constant for the Hardy inequality as we already
know before.
\end{proof}

The same proof could be adapted with minor requirements to prove
the following

\begin{theorem}
Let $u$ be a function in $C_c^{\infty}(\R^{n} \setminus \{0\})$
and $p > n > 1$. Then,
\begin{equation}
\label{HR}
\int_{\R^{n}}\frac{\left|u(x)\right|^{p}}{\left\|x\right\|^{p}}dx
\leq \left(\frac{p}{p-n}\right)^{p}\int_{\R^{n}}\left\|\nabla
u(x)\right\|^{p}dx,
\end{equation}
where $\left(\frac{p}{p-n}\right)^{p}$ is the sharp constant.
\end{theorem}

\subsection{The case $(p>n, n=1)$}

Now we are going to prove a sharp version of the Hardy's
inequality in one dimension.

\begin{theorem}
Let $u$ be a nonnegative function in $C_c^{\infty}(\R_+)$, which
is non-identically zero, $p>1$ and set
$$
  \eta(x):= \int_{0}^{x} u(t)\, dt, \quad \quad \Big( \frac{d\eta(x)}{dx} \equiv \eta'(x) = u(x) \Big).
$$
Then, we have
$$
  \int_{0}^{\infty} \frac{\eta(x)^p}{x^p} \, dx
  < \left(\frac{p}{p-1}\right)^{p}\int_{0}^{\infty} u(x)^{p} \,
  dx,
$$
where $\left(\frac{p}{p \,-1}\right)^{p}$ is the sharp constant.
\end{theorem}

\begin{proof}

First, we observe that $\eta(0)=0$. The proof follows almost the
same lines as before. Indeed, we have
$$
\begin{aligned}
  \int_{0}^{\infty} \eta(x)^p \, \frac{1}{x^p} \, dx&=
  \frac{1}{1-p}\int_{0}^{\infty} \eta(x)^{p} \, (x^{1-p})' \, dx \\[5pt]
&=
  \frac{p}{p-1}\int_{0}^{\infty}\left(\frac{\eta(x)}{x}\right)^{p-1} \, \eta'(x) \,
  dx,
\end{aligned}
$$
where we have integrated by parts and used that
$$
  \begin{aligned}
\lim_{x\rightarrow 0} \eta(x)^{p} \; x^{1-p}&=0, \qquad (\eta(0)=0), \\[5pt]
\lim_{x\rightarrow \infty} \eta(x)^{p} \; x^{1-p}&= 0, \qquad
(p>1).
\end{aligned}
$$
Therefore, we obtain
\begin{equation} \int_{0}^{\infty}\left(\frac{\eta(x)}{x}\right)^{p} \, dx=
\frac{p}{p-1}\int_{0}^{\infty}\left(\frac{\eta(x)}{x}\right)^{p-1}
\, u(x) \, dx.
\end{equation}
Now applying the Young's inequality, it follows that
$$
  \int_{0}^{\infty}\left(\frac{\eta(x)}{x}\right)^{p} \, dx
  <
  \lambda^{q} \int_{0}^{\infty}\left(\frac{\eta(x)}{x}\right)^{p} \, dx
  +\frac{1}{(p-1)\lambda^{p}}\int_{0}^{\infty} u(x)^{p} \, dx,
$$
i.e.,
$$
  \int_{0}^{\infty}\left(\frac{\eta(x)}{x}\right)^{p} \, dx
  < \frac{1}{(p-1)\lambda^{p}(1-\lambda^{q})}\int_{0}^{\infty}
  u(x)^{p}dx.
$$
One remarks that, the equality happens in the Young's inequality
before if, and only if
$$
  \left(\frac{\eta(x)}{x}\right)^{(p-1) q}= u(x)^{p}.
$$
Therefore, it follows that $u$ must be a positive constant, which
is a contradiction since $u$ has compact support and is
non-identically zero.

\medskip
Finally, we define for $p$ and thus $q$ fixed,
$$
  f(\lambda)= \frac{1}{(p-1)\lambda^{p}(1-\lambda^{q})}.
$$
Then, we find that the minimum point of $f$ is
$\lambda_{0}=q^{-q}$, and moreover the minimal value is
$$
  f(\lambda_{0})=\left(\frac{p}{p-1}\right)^{p}.
$$
Consequently, we obtain the sharp Hardy's inequality, that is
$$
  \int_{0}^{\infty}\left(\frac{\eta(x)}{x}\right)^{p} \, dx
  <\left(\frac{p}{p-1}\right)^{p}\int_{0}^{\infty} u(x)^{p} \, dx.
$$
\end{proof}

\section{The Caffarelli-Kohn-Nirenberg inequality} \label{CKNI}

In order to show the inequality due to Caffarelli-Kohn-Nirenberg,
we consider first two lemmas, which are particular cases. The
former is this inequality for $b=a+1$, the second one is for
$b=a$. The general result follows by an interpolation argument.

\begin{lemma}
\label{L1} Let $u$ be a function in $C_{c}^{\infty}(\R^{n})$,
$n\geq 3$,
$$
  -\infty<a<\frac{n-2}{2} \quad \text{and} \quad b=a+1.
$$
Then, we have
\begin{equation}
\label{CKNA1}
\int_{\R^{n}}\frac{\left|u(x)\right|^{2}}{\left\|x\right\|^{2b}}
\, dx \leq
C_{a+1} \int_{\R^{n}}\frac{\left\|\nabla
u(x)\right\|^{2}}{\left\|x\right\|^{2a}} \, dx,
\end{equation}
where
$$
  C_{a+1}= \frac{4}{(n-2-2a)^{2}}
$$
is the sharp constant.
\end{lemma}

\begin{proof}
Let $W:\R^{n} \setminus \{0\} \to \R^{n}$ be a smooth vector value
function defined as
$$
  W(x):=\left(\frac{1}{2b-n}\right)\frac{x}{\left\|x\right\|^{2b}}.
$$
This vector field is well defined since $n\neq 2b$. Indeed, if $b=
n/2$, then we must have $a= (n-2)/2$, which contradicts the
hypothesis. Moreover, a straightforward calculation shows that
$$
  \dive W(x)=-\frac{1}{\left\|x\right\|^{2b}}.
$$
Now, the integral of the left side of $(\ref{CKNA1})$ can be
written of the following way
$$
  \int_{\R^{n}}\frac{\left|u(x)\right|^{2}}{\left\|x\right\|^{2b}} \, dx=
  -\int_{\R^{n}}\left|u(x)\right|^{2} \dive W(x) \, dx.
$$
Then, applying the Gauss-Green Theorem and using the Young's
inequality, we obtain
\begin{equation}
\label{CKNAA}
  \begin{aligned}
  \int_{\R^{n}} \frac{\left|u(x)\right|^{2}}{\left\|x\right\|^{2b}} \, dx&=
  2 \int_{\R^{n}} \left( \frac{\left|u(x)\right| \,
  W(x)}{\left\|x\right\|^{-a}}\right) \cdot \frac{ \left(\nabla \left|u(x)\right|\right)}
  {\left\|x\right\|^{a}} \\[5pt]
  &\leq \alpha^{2}\int_{\R^{n}}\frac{\left|u(x)\right|^{2}\left\|W(x)\right\|^{2}}{\left\|x\right\|^{-2a}} \, dx
  + \alpha^{-2} \int_{\R^{n}}\frac{\left\|\nabla u(x)\right\|^{2}}{\left\|x\right\|^{2a}} \,
  dx,
  \end{aligned}
\end{equation}
where $\alpha$ is a positive real number. Now, we observe that
$b=a+1$ implies
$$
  \frac{\left\|W(x)\right\|^{2}}{\left\|x\right\|^{-2a}}
  =\frac{\|x\|^2}{(2b-n)^{2}\left\|x\right\|^{4a+4}} \;
  \frac{1}{\|x\|^{-2a}}
  =\frac{1}{(n-2-2a)^{2}\left\|x\right\|^{2b}}.
$$
Therefore, we have from \eqref{CKNAA}
$$
  \left(1-\frac{\alpha^{2}}{(n-2-2a)^{2}}\right)
  \int_{\R^{n}}\frac{\left|u(x)\right|^{2}}{\left\|x\right\|^{2b}}dx
  \leq \alpha^{-2}\int_{\R^{n}}\frac{\left\|\nabla u(x)\right\|^{2}}{\left\|x\right\|^{2a}} \, dx
$$
and so
$$\int_{\R^{n}}\frac{\left|u(x)\right|^{2}}{\left\|x\right\|^{2b}}dx
\leq f(\alpha;n,a) \int_{\R^{n}}\frac{\left\|\nabla
u(x)\right\|^{2}}{\left\|x\right\|^{2a}}dx,
$$
where
\begin{equation}
  f(\alpha;n,a)= \frac{(n-2-2a)^{2}}{\alpha^{2}(n-2-2a)^{2}-\alpha^{4}}.
\end{equation}

Analogously, we set $\kappa= (n-2-2a)^{2}$ and look for the
critical points of $f(\alpha)$. Thus we proceed as before and
obtain the critical point
$$
  \alpha_{0}= \sqrt{\frac{\kappa}{2}}
$$
and also $f''(\alpha_0)>0$. Moreover, we have
$$
  f(\sqrt{\kappa/2})= \frac{4}{(n-2-2a)^{2}},
$$
which is the best constant for this inequality.
\end{proof}

\bigskip
Now we are going to prove a second lemma, which is the
Caffarelli-Kohn-Nirenberg's inequality when $b=a$. In this case,
we do not follow the same ideas as before, moreover we do not
recover the best constant.

\begin{lemma}
\label{L2} Let $u$ be a function in $C_{c}^{\infty}(\R^{n})$, $n
\geq 3$,
$$
  -\infty < 2 \, a < n-2, \quad b=a \quad \text{and} \quad p=\frac{2n}{n-2}=2^*.
$$
Then, we have
\begin{equation}
\label{CKN}
   \left(\int_{\R^{n}}\frac{\left|u(x)\right|^{p}}{\left\|x\right\|^{a p}}
   \, dx \right)^{\frac{2}{p}}
   \leq
   C_{a^\pm}
   \int_{\R^{n}}\frac{\left\|\nabla
   u(x)\right\|^{2}}{\left\|x\right\|^{2a}} \, dx,
\end{equation}
where
$$
  C_{a^+}= K(n,p)^2 \, \Big( \frac{n-2}{n - 2 - 2a} \Big)^{2}, \quad
  C_{a^-}= K(n,p)^2 \, \Big( \frac{n-2-4a}{n - 2 - 2a} \Big)^{2},
$$
and $a^+$ stands for $a\geq 0$ similarly $a^-$ for $a \leq 0$.
\end{lemma}

\begin{proof}
We begin applying the Sobolev inequality \eqref{SI} for $f(x)=
u(x) / \|x\|^{a}$ and without misunderstanding $(p= 2)$.
Therefore, we obtain
\begin{equation}
\label{CKN2}
   \left(\int_{\R^{n}} \frac{|u(x)|^p}{\left\|x\right\|^{ap}} \, dx \right)^{\frac{2}{p}}
   \leq K(n,p)^{2}\int_{\R^{n}}\left\|\nabla\left(\frac{u(x)}{\left\|x\right\|^{a}}\right)\right\|^{2}dx.
\end{equation}
Now, we analyze the right-hand side of the above inequality.
First, we observe that

%
$$
  \|\nabla \Big( \frac{u(x)}{\|x\|^a} \Big)\|^2=
  \frac{\| \nabla u(x)\|^2}{\|x\|^{2a}}
  + a^2 \frac{|u(x)|^2}{\|x\|^{2(a+1)}}
  - \frac{2 \, a \, u(x)}{\|x\|^{2a+2}} \nabla u(x) \cdot x,
$$
and thus
\begin{equation}
\label{CKNL2I}
  \begin{aligned}
  \int_{\R^n} \|\nabla \Big( \frac{u(x)}{\|x\|^a} \Big)\|^2 \, dx &=
  \int_{\R^n} \frac{\| \nabla u(x)\|^2}{\|x\|^{2a}} \, dx \\[5pt]
  &+  \int_{\R^n} a^2 \frac{|u(x)|^2}{\|x\|^{2(a+1)}} \, dx
  - \int_{\R^n} \frac{2 \, a \, u(x)}{\|x\|^{2a+2}} \nabla u(x) \cdot x \, dx
  \\[5pt]
  &= I_1 + I_2 + I_3,
  \end{aligned}
\end{equation}
with the obvious notations. For $I_1$ term there is nothing to do,
and for $I_2$ we apply Lemma \ref{L1}, then we have
\begin{equation}
\label{CKNL2I2}
   I_2 \leq a^2 \, C_{a+1} \int_{\R^n} \frac{\| \nabla u(x)\|^2}{\|x\|^{2a}} \,
   dx.
\end{equation}
It remains to consider the $I_3$ term. We divide in two cases:
$a<0$ and $a > 0$. For the former, we have
\begin{equation}
\label{CKNL2I22}
  \begin{aligned}
  I_3 &\leq - a \lambda^2 \int_{\R^n}
  \frac{|u(x)|^2}{\|x\|^{2(a+1)}} \, dx
  \, - \,  a \lambda^{-2} \int_{\R^n}
  \frac{\| \nabla u(x)\|^2}{\|x\|^{2a}} \, dx \\[5pt]
  &\leq -a \, \big(\lambda^2 \, C_{a+1} + \lambda^{-2} \big)
  \int_{\R^n} \frac{\| \nabla u(x)\|^2}{\|x\|^{2a}} \, dx,
  \end{aligned}
\end{equation}
where we have used Young's inequality \eqref{YI} and Lemma
\ref{L1}. Therefore, from \eqref{CKN2}--\eqref{CKNL2I22}, we
obtain
$$
  \left(\int_{\R^{n}}\frac{\left|u(x)\right|^{p}}{\left\|x\right\|^{a p}}
   \, dx \right)^{\frac{2}{p}}
   \leq f_{-}(\lambda;n,a) \,
  \int_{\R^{n}}\frac{\left\|\nabla
   u(x)\right\|^{2}}{\left\|x\right\|^{2a}} \, dx,
$$
where
$$
  f_-(\lambda;n,a)= K(n,p)^2 \Big( 1 + a^2 C_{a+1} - a
  \big(\lambda^2 C_{a+1} + \lambda^{-2} \big) \Big).
$$
Analogously, we proceed for the second case when $a > 0$, and
obtain
$$
  \left(\int_{\R^{n}}\frac{\left|u(x)\right|^{p}}{\left\|x\right\|^{a p}}
   \, dx \right)^{\frac{2}{p}}
   \leq f_{+}(\lambda;n,a) \,
  \int_{\R^{n}}\frac{\left\|\nabla
   u(x)\right\|^{2}}{\left\|x\right\|^{2a}} \, dx,
$$
where
$$
  f_+(\lambda;n,a)= K(n,p)^2 \Big( 1 + a^2 C_{a+1} + a
  \big(\lambda^2 C_{a+1} + \lambda^{-2} \big) \Big).
$$

\bigskip
Finally, we proceed to minimize $f_\pm$ with respect to $\lambda$.
Thus we obtain the minimal point
$$
  \lambda_0= \Big(\frac{1}{C_{a+1}} \Big)^{1/4},
$$
such that $f'_\pm(\lambda_0)= 0$ and $f''_\pm(\lambda_0)>0$.
Moreover, it follows that
$$
  f_\pm(\lambda_0)= K(n,p)^2 \Big( 1 \pm a \sqrt{C_{a+1}} \Big)^2.
$$
Consequently, we obtain
$$
  \left(\int_{\R^{n}}\frac{\left|u(x)\right|^{p}}{\left\|x\right\|^{ap}}dx\right)^{\frac{2}{p}}
  \leq
  K(n,p)^{2} \left(1 \pm a \sqrt{C_{a+1}} \right)^{2}\int_{\R^{n}}\frac{\left\|\nabla
  u(x)\right\|^{2}}{\left\|x\right\|^{2a}}dx.
$$
\end{proof}

\bigskip
We finish this section with the proof of
Caffarelli-Kohn-Nirenberg's inequality.

\begin{theorem}
\label{CKNT} Let $u$ be a function in $C_c^{\infty}(\R^{n})$ for
$n\geq 3$ and assume that
$$
  -\infty \leq a <\frac{n-2}{2}, \quad a \leq b \leq a+1
  \quad \text{and} \quad
  p=\frac{2n}{n-2 + 2(b-a)}.
$$
Then, we have for each $\theta \in [0,1]$
\begin{equation}
\label{CKN4}
\left(\int_{\R^{n}}\frac{\left|u(x)\right|^{p}}{\left\|x\right\|^{bp}}
\,dx\right)^{\frac{2}{p}} \leq C_{a^\pm}^\alpha \, C_{a+1}^\beta
\, \int_{\R^{n}}\frac{\left\|\nabla
u(x)\right\|^{2}}{\left\|x\right\|^{2a}} \, dx,
\end{equation}
where $C_{a^\pm}$ and $C_{a+1}$ are respectively the constants for
$b=a$ and $b=a+1$, and
$$
  \alpha= \frac{2n \, \theta}{(n-2)p}, \quad \beta=
  \frac{2\, (1-\theta)}{p}.
$$
\end{theorem}

\begin{proof}
We are going to obtain \eqref{CKN4} by interpolation between the
two previous results obtained in Lemma \ref{L1} and Lemma
\ref{L2}. First, we write for each $0<\theta<1$ the exponent $p$
as
\begin{equation}
\label{EXP1}
  p= 2(1-\theta) + 2^{*}\theta.
\end{equation}
We recall that $2^{*}$ is the Sobolev conjugate of 2, i.e. $2^*=
2n/(n-2)$. Since we have
$$
  p= \frac{2n}{n-2+2(b-a)},
$$
it follows from \eqref{EXP1} after a straightforward algebraic
calculation that
\begin{equation}
\label{EXP2} b= a+1-\left(\frac{n\theta}{n-2+2\theta}\right)
\end{equation}
and also
\begin{equation}
\label{EXP3} b \,p= 2 \, (1-\theta)(a+1) + 2^{*} \, \theta \, a.
\end{equation}
Therefore, from \eqref{EXP1}--\eqref{EXP3} we could write
$$
\begin{aligned}
\int_{\R^{n}}\frac{|u(x)|^{p}}{\left\|x\right\|^{bp}} \, dx &=
\int_{\R^{n}}\frac{|u(x)|^{2(1-\theta)+2^{*}\theta}}{\left\|x\right\|^{2(1-\theta)(a+1)+2^{*}
\theta a}} \, dx
\\[5pt]
&\leq
\left(\int_{\R^{n}}\frac{|u(x)|^{2}}{\left\|x\right\|^{2(a+1)}} \,
dx \right)^{1-\theta}
\left(\int_{\R^{n}}\frac{|u(x)|^{2^{*}}}{\left\|x\right\|^{2^{*}a}}
\, dx\right)^{\theta},
\end{aligned}
$$
where we have applied Hölder's inequality with
$$
  \tilde{p}= \frac{1}{1-\theta}
  \quad \text{and} \quad
  \tilde{q}= \frac{1}{\theta}.
$$

The proof follows from the Lemma \ref{L1} and Lemma \ref{L2}, that
is

$$\left(\int_{\R^{n}}\frac{\left|u(x)\right|^{p}}{\left\|x\right\|^{bp}}dx\right)^{\frac{2}{p}}\leq
C_{a^\pm}^{\frac{2^{*}\theta}{p}}C_{a+1}^{\frac{2(1-\theta)}{p}}\int_{\R^{n}}\frac{\left\|\nabla
u(x)\right\|^{2}}{\left\|x\right\|^{2a}}dx.$$
\end{proof}

\begin{remark}
In the previous theorem, instead of $C_{a^\pm}$ and $C_{a+1}$ we
could consider respectively $C(n,b,2^*)$ and $C(n,b,2)$. Thus for
each $\theta \in [0,1]$ fixed, we obtain the constant of the
Caffarelli-Kohn-Nirenberg's inequality, see \cite{CATRINA}, that
is
$$
  C(n,b,2^*)^\alpha \, C(n,b,2)^\beta
$$
with $\alpha$ and $\beta$ given by Theorem \ref{CKNT}
\end{remark}

\section{The Rellich inequality} \label{RIS}

In this last section we prove Rellich's inequality, which is a
second order generalization of the Hardy's inequality.

\begin{theorem}
Let $u$ be a function in $C_{c}^{\infty}(\R^{n})$ and $n>4$. Then,
we have
\begin{equation}
\label{RR}
    \int_{\R^{n}}\frac{\left|u(x)\right|^{2}}{\left\|x\right\|^{4}} \,
    dx
    \leq
    \left(\frac{4}{n(n-4)}\right)^2 \int_{\R^{n}}\left|\Delta
u(x)\right|^{2}dx,
\end{equation}
where $\left(\frac{4}{n(n-4)}\right)^2$ is the sharp constant.

\end{theorem}

\begin{proof}
The proof has 3 parts:

1. First, we use the particular case of Caffarelli-Kohn-Nirenberg
inequality for $b=a+1$, with $a=1$, so
\begin{equation}
\label{RI1}
\int_{\R^{n}}\frac{\left|u(x)\right|^{2}}{\left\|x\right\|^{4}}dx\leq
M\int_{\R^{n}} \frac{\left\|\nabla
u\right\|^{2}}{\left\|x\right\|^{2}} \, dx,
\end{equation}
where
$$
  M:=\frac{4}{(n-4)^{2}}.
$$

2. Now, we observe the right side of \eqref{RI1}, and we can note
that

$$
  M \int_{\R^{n}}\frac{\left\|\nabla u\right\|^{2}}{\left\|x\right\|^{2}}
  =
  M \int_{\R^{n}} \frac{2 \, u}{\left\|x\right\|^{4}} \; x \cdot \nabla u \, dx
  -M \int_{\R^{n}}\frac{u}{\left\|x\right\|^{2}} \, \Delta u \, dx.
$$

Thus using the Young's inequality, we obtain

$$
\begin{aligned}
  M\int_{\R^{n}}\frac{\left\|\nabla u\right\|^{2}}{\left\|x\right\|^{2}}
  &\leq
  \frac{M^{2}}{2\lambda^{2}}\int_{\R^{n}}\frac{\left\|\nabla u\right\|^{2}}{\left\|x\right\|^{2}}
  +\frac{M\lambda^{2}}{2}\int_{\R^{n}}\left|\Delta u\right|^{2}
  \\[5pt]
  &+\frac{M^{2}}{\mu^{2}}\int_{\R^{n}}\frac{\left\|\nabla
  u\right\|^{2}}{\left\|x\right\|^{2}}+\mu^{2}M\int_{\R^{n}}\frac{\left\|\nabla
  u\right\|^{2}}{\left\|x\right\|^{2}},
\end{aligned}
$$

where $\lambda$, $\mu>0$. So, we have
\begin{equation}
\label{RI2}
    \int_{\R^{n}}\frac{\left|u(x)\right|^{2}}{\left\|x\right\|^{4}} \,dx
     \leq P(\lambda,\mu;M)\int_{\R^{n}}\left|\Delta u(x)\right|^{2} \,
     dx,
\end{equation}
where
$$
  P(\lambda,\mu;M):= M \, \frac{\lambda}{2} \,
  \left(1-\frac{M}{2\lambda^{2}}-\frac{M}{\mu^{2}}-\mu^{2}\right)^{-1}.
$$

3. Finally, we can proceed to obtain the critical point
$(\lambda_{0},\mu_{0})$ of $P$. In fact, replacing this critical
point in the determinant of the Hessian of $P$, i.e. matrix $H$,
we have
$$
   \det H(\lambda_{0},\mu_{0}) \geq 0.
$$
Moreover, we obtain
$$
  P(\lambda_{0},\mu_{0})=\frac{16}{n^{2}(n-4)^{2}}.
$$
Consequently, from \eqref{RI2}, it follows that

$$\int_{\R^{n}}\frac{\left|u(x)\right|^{2}}{\left\|x\right\|^{4}}dx
\leq \frac{16}{n^{2}(n-4)^{2}}\int_{\R^{n}}\left|\Delta
u(x)\right|^{2}dx.$$
\end{proof}

\section*{Acknowledgements}
The first author is supported by FAPERJ by the grant 2009.2848.0.
The second author were partially supported by FAPERJ through the
grant E-26/ 111.564/2008 entitled {\sl ``Analysis, Geometry and
Applications''}.


%

\end{document}